\documentclass[runningheads]{llncs}
\usepackage{amssymb}
\usepackage{graphicx}
\usepackage{url}
\usepackage{amsmath} 
\usepackage{mathrsfs, multirow, subfigure}
\usepackage{ifthen} 
\usepackage{amsfonts}
\usepackage[usenames]{color}
\usepackage[normalem]{ulem}
\usepackage{tcolorbox}
\usepackage{soul}
\usepackage[ruled,vlined]{algorithm2e}
\usepackage{float}

\newcommand{\keywords}[1]{\par\addvspace\baselineskip
\noindent\keywordname\enspace\ignorespaces#1}

\usepackage[sort&compress]{natbib}
\bibpunct{(}{)}{;}{a}{}{,} 
\makeatletter
\renewcommand\bibsection%
{
  \section*{\refname
    \@mkboth{\MakeUppercase{\refname}}{\MakeUppercase{\refname}}}
}
\makeatother

\newtheorem{rcond}{Regularity Conditions}

\newcommand{\prob}{\mathsf{P}}

\newcommand{\nm}{{\sf N}}

\newcommand{\RR}{\mathbb{R}}

\newcommand{\TT}{\mathbb{T}}

\newcommand{\ZZ}{\mathbb{Z}}

\newcommand{\lPi}{\underline{\Pi}}
\newcommand{\uPi}{\overline{\Pi}}

\newcommand{\uGamma}{\overline{\Gamma}}

\begin{document}
\title{Large-sample theory for inferential models: a possibilistic Bernstein--von Mises theorem}
\titlerunning{A possibilistic Bernstein--von Mises theorem}
\author{Ryan Martin\and 
Jonathan P. Williams}
\authorrunning{R.~Martin and J.~P.~Williams}
\institute{North Carolina State University \\ Department of Statistics \\ Raleigh, North Carolina 27695, U.S.A. \\ 
\email{\{rgmarti3, jwilli27\}@ncsu.edu}\\
}
\maketitle              
\begin{abstract}
The inferential model (IM) framework offers alternatives to the familiar probabilistic (e.g., Bayesian and fiducial) uncertainty quantification in statistical inference.  Allowing this uncertainty quantification to be imprecise makes it possible to achieve exact validity and reliability.  But is imprecision and exact validity compatible with attainment of the classical notions of statistical efficiency?  The present paper offers an affirmative answer to this question via a new possibilistic Bernstein--von Mises theorem that parallels a fundamental result in Bayesian inference.  Among other things, our result demonstrates that the IM solution is asymptotically efficient in the sense that its asymptotic credal set is the smallest that contains the Gaussian distribution whose variance agrees with the Cram\'er--Rao lower bound.  

\keywords{Asymptotics; Bayesian; belief; fiducial; relative likelihood}
\end{abstract}

\section{Introduction}
\label{S:intro}

\citet{efron.cd.discuss} writes that ``the most important unresolved problem in statistical inference is the use of Bayes theorem in the absence of prior information.''  Numerous attempts have been made, including Bayesian inference with default priors, from Laplace's flat priors to the refinements by \citet{jeffreys1946} and \citet{bergerbernardosun2009}, fiducial inference in Fisher's original sense \citep{fisher1935a, zabell1992} and its generalizations \citep{fraser1968, hannig.review, dempster1967, dempster2008}, and the various imprecise-probabilistic proposals \citep[e.g.,][]{berger1984, walley1991, dubois2006, augustin.etal.bookchapter}.  One thing that many of these solutions offer is a large-sample result which states, roughly, that the associated credible sets achieve the nominal frequentist coverage probability when the sample size, $n$, is large.  Statisticians insist on methods that are demonstrably reliable in the sense of tending to report ``correct inferences'' across repeated uses, so results like this are high-value.  The {\em r\'esultat extraordinaire} along these lines is the {\em Bernstein--von Mises theorem} which states the following: under certain conditions, the Bayesian or fiducial posterior distribution is approximately Gaussian, centered at the maximum likelihood estimator, and is efficient in the sense that its variance equals the Cram\'er--Rao lower bound; here ``approximately'' is in the sense that the total variation distance between the two distributions is vanishing in probability as $n \to \infty$.  For details, see \citet[][Ch.~10]{vaart1998}, \citet[][Ch.~4]{ghosh-etal-book}, or \citet{hannig.review}. 

Still, none of the Bayesian- or fiducial-like solutions mentioned above have proved to be fully satisfactory.  One possible explanation is that large-sample confidence sets don't tell the whole story.  Indeed, the false confidence theorem \citep{balch.martin.ferson.2017, martin.nonadditive} implies that, at every $n$, there exists false hypotheses about the unknown to which the (Bayesian or fiducial) posterior will tend to assign relatively high probability---that is, the posterior probability assignments have an inherent unreliability.  To avoid this pitfall, it's necessary that the Bayesian or fiducial precise probability be replaced by a data-dependent imprecise probability.  The inferential model (IM) framework, including the original developments in \citet{imbasics, imbook} and the more recent generalizations in \citet{imchar, martin.partial, martin.partial2}, does just this.  Moreover, the IM is provably valid in a sense that implies it's safe from false confidence and, e.g., it provides exact finite-sample confidence sets.  While the IM solution's efficiency compared to existing solutions---Bayesian, frequentist, or otherwise---has been determined empirically in specific applications, a general theoretical statement concerning its efficiency is lacking.  In this paper, for the class of IMs in \citet{plausfn, gim, martin.partial2}, we establish a possibilistic version of the Bernstein--von Mises theorem that parallels those mentioned above.  Among other things, this result confirms our conjecture that, despite the genuine imprecision in the IM's output needed to ensure exact validity, there's no loss of efficiency asymptotically.

\section{Background}

\subsection{Possibility theory}
\label{SS:possibility}

Possibility theory is among the simplest imprecise probability theories, corresponding to {\em consonant} belief structures \citep[e.g.,][Ch.~10]{shafer1976}.  Other key references include \citet{zadeh1978}, \citet{dubois.prade.book}, and \citet{dubois2006}. 
The simplicity of possibilistic uncertainty quantification comes from its parallels to precise probability theory.  A necessity--possibility measure pair $(\lPi, \uPi)$ that's intended to quantify uncertainty about an uncertain $Z$ in $\ZZ$ is determined by a possibility contour $\pi: \ZZ \to [0,1]$, with $\sup_{z \in \ZZ} \pi(z) = 1$, via the rules
\[ \uPi(B) = \sup_{z \in B} \pi(z) \quad \text{and} \quad \lPi(B) = 1 - \sup_{z \not\in B} \pi(z), \quad B \subseteq \ZZ. \]
So, where ordinary probability is determined by integrating a density function, possibility is determined by optimizing a contour.  The above values are often interpreted subjectively as (coherent) upper and lower probabilities associated with the proposition ``$Z \in B$'' or as Shaferian degrees of belief.

One way to elicit a possibility measure, specifically relevant to us here, is via the {\em probability-to-possibility transform} \citep[e.g.,][]{dubois.etal.2004, hose2022thesis}.  If $p$ is a probability density function, which determines a random variable $Z \sim \prob$, then the probability-to-possibility transform defines $\pi$ as 
\[ \pi(z) = \prob\{ p(Z) \leq p(z) \}, \quad z \in \ZZ. \]
This defines the ``best approximation'' 
of $\prob$ by a possibility measure in the sense that its corresponding credal set is the smallest one that contains $\prob$.  If $p$ is the (multivariate) standard normal density on $\ZZ$, the probability-to-possibility transform has contour defined as $\gamma$ as shown in Figure~\ref{fig:bernoulli}(a); the corresponding Gaussian possibility measure defined via optimization is denoted by $\uGamma$.

\subsection{Inferential models}
\label{SS:im}

The original IM constructions \citep[e.g.,][]{imbasics, imbook} were formulated in terms of (nested) random sets and, hence, the connection to possibility theory was indirect.  A more recent IM construction introduced in \citet{martin.partial2}---see, also, \citet{plausfn, gim}---directly defines the IM's possibility contour using a version of the probability-to-possibility transform.  An advantage of this alternative perspective and construction is that it avoids the ambiguity in specifying both a data-generating equation and a so-called ``predictive random set.''  



The statistical model assumes that $X^n=(X_1,\ldots,X_n)$ consists of iid samples from a distribution $\prob_\Theta$ depending on an unknown/uncertain $\Theta \in \TT$ to be inferred.  The model and observed data $X^n=x^n$ together determine a likelihood function $\theta \mapsto L_{x^n}(\theta)$ and a corresponding relative likelihood function
\[ 
R(x^n,\theta) = \frac{L_{x^n}(\theta)}{\sup_\vartheta L_{x^n}(\vartheta)}.
\]
The relative likelihood itself defines a data-dependent possibility contour that has been widely studied \citep[e.g.,][]{shafer1982, wasserman1990b, denoeux2006, denoeux2014}.  Most appealing about the likelihood-based possibility contour is its shape: peak at the maximum likelihood estimator $\hat\theta_{x^n}$ and consistent with Fisher's suggested likelihood-based preference order on the parameter space.  But it lacks a standard scale, i.e., what constitutes a ``small'' relative likelihood depends on aspects of the individual application.  A (literally) uniform scale of interpretation across applications can easily be obtained via what \citet{martin.partial} calls ``validification,'' which is a sort of possibilistic transform.  In particular, for observed data $X^n=x^n$, the possibilistic IM's contour is
\begin{equation}
\label{eq:contour}
\pi_{x^n}(\theta) = \prob_\theta\{ R(X^n,\theta) \leq R(x^n, \theta) \}, \quad \theta \in \TT,
\end{equation}
and the corresponding possibility measure is 
\[ \uPi_{x^n}(H) = \sup_{\theta \in H} \pi_{x^n}(\theta), \quad H \subseteq \TT. \]
Critical to the IM developments is the so-called {\em validity property} which can be succinctly described by the following expression:
\begin{equation}
\label{eq:valid}
\sup_{\Theta \in \TT} \prob_\Theta\bigl\{ \pi_{X^n}(\Theta) \leq \alpha \bigr\} \leq \alpha, \quad \text{for all $\alpha \in [0,1]$}. 
\end{equation}
This is precisely the universal scaling that the relative likelihood itself is missing, i.e., ``$\pi_{x^n}(\theta) \leq \alpha$'' has the same inferential meaning/force in every application.  It also implies that the IM is safe from false confidence.  Finally, it has some important and familiar statistical consequences:
\begin{itemize}
\item a test that rejects the hypothesis ``$\Theta \in H$'' when $\uPi_{x^n}(H) \leq \alpha$ will control the frequentist Type~I error rate at level $\alpha$, and 
\item the set $C_\alpha(x^n) = \{\theta \in \TT: \pi_{x^n}(\theta) > \alpha\}$ is a $100(1-\alpha)$\% frequentist confidence set in the sense that its coverage probability is at least $1-\alpha$.
\end{itemize}
Given that the IM is inherently imprecise and offers finite-sample guarantees, the reader might think that this comes with an associated loss of efficiency compared to those methods whose justification relies on asymptotic considerations.  The result in the next section aims to debunk this myth.

\section{A possibilistic Bernstein--von Mises theorem}

\subsection{Preview}

To build some intuition, consider a simple exactly Gaussian model where $X^n = (X_1,\ldots,X_n)$ is an iid sample from $\nm(\Theta,\sigma^2)$, where the mean $\Theta$ is unknown  but $\sigma^2$ is known.  The maximum likelihood estimator is $\hat\theta_{x^n}=\bar x$, the sample mean, and it is 
not difficult to show that the IM's possibility contour is 
\[ \pi_{x^n}(\theta) = \gamma\{ n^{1/2}(\theta-\hat\theta_{x^n}) / \sigma \}, \quad \theta \in \RR, \]
where $\gamma$ is the Gaussian possibility contour in Section~\ref{SS:possibility}.  Alternatively, if we switch to a local parametrization of $\theta$, i.e., $\theta = \hat\theta_{x^n} + \sigma n^{-1/2} z$, then 
\[ \check\pi_{x^n}(z) := \pi_{x^n}(\hat\theta_{x^n} + \sigma n^{-1/2} z) = \gamma(z), \quad z \in \RR. \]
For a generic $H \subseteq \TT$, the corresponding possibility measure is 
\begin{align*}
\uPi_{x^n}(H) & = 
\uGamma\bigl\{ n^{1/2}(H - \hat\theta_{x^n})/\sigma \bigr\}, 
\end{align*}
where $\uGamma$ is the Gaussian possibility measure from Section~\ref{SS:possibility} and, for constants $a$ and $b$, the set $aH+b$ is defined as $aH + b = \{a\theta + b: \theta \in H\}$. 

Once we step away from an exact Gaussian model, there is no longer an exact correspondence between the possibilistic IM's solution and the Gaussian possibility measure.  It's a similar story in the Bayesian and (generalized) fiducial case.  But the probabilistic Bernstein--von Mises theorem implies that, under certain conditions, as $n \to \infty$, the suitably centered and scaled Bayesian posterior distribution will be approximately Gaussian.  Our main result shows that the same holds true for the possibilistic IM solution: under certain regularity conditions, a suitably centered and scaled version of the IM's possibility measure will converge in a strong sense to the Gaussian possibility measure.

\subsection{Main result}
\label{SS:main}

The regularity conditions stated below are exactly those assumed in textbook treatments of the asymptotic normality of maximum likelihood estimators---in particular, \citet[][Theorem~7.63]{schervish1995}. 

\begin{rcond}
The probability measures $\prob_\theta$, indexed by $\theta \in \TT \subseteq \RR^d$, have density functions $p_\theta$ with respect to a fixed dominating measure.  Recall that $\Theta$ denotes the special ``true'' parameter value.  
\begin{enumerate}
\item The support of $\prob_\theta$ doesn't depend on $\theta$.
\item $\theta \mapsto \log p_\theta(x)$ is twice continuously differentiable for almost all $x$.
\item Differentiation can be passed under the integral sign. 
\item The second derivative of $\theta \mapsto \log p_\theta(X_1)$ is Lipschitz continuous in a neighborhood of $\Theta$ with Lipschitz constant that has finite $\prob_\Theta$-expectation.
\end{enumerate}
\end{rcond}


Under the above regularity conditions, we establish a possibilistic Bernstein--von Mises theorem for the IM solution presented in Section~\ref{SS:im}.  That is, a centered and scaled version of the IM's possibility contour, namely, 
\begin{equation}
\label{eq:im.scaled}
\check\pi_{X^n}(z) = \pi_{X^n}\bigl( \hat\theta_{X^n} + J_{X^n}^{-1/2} z \bigr), \quad z \in \RR^d, 
\end{equation}
where $J_{X^n}$ denotes the observed Fisher information matrix, can be well approximated by the Gaussian possibility contour.  

\begin{theorem}
\label{thm:bvn}
Let $X^n = (X_1,\ldots,X_n)$ consist of iid observations from $\prob_\Theta$, where $\{\prob_\theta: \theta \in \TT\}$ is the posited model and $\Theta$ is in the interior of $\TT$.  If the above regularity conditions hold at $\Theta$, and the maximum likelihood estimator $\hat\theta_{X^n}$ is consistent, then the possibilistic IM's solution satisfies 
\begin{equation}
\label{eq:contour.limit}
\sup_{z \in K} \bigl| \check\pi_{X^n}(z) - \gamma(z) \bigr| \to 0, \quad \text{in $\prob_\Theta$-probability as $n \to \infty$}, 
\end{equation}
where $K \subset \RR^d \setminus \{0\}$ is an arbitrary compact set, $\check\pi_{X^n}$ is as in \eqref{eq:im.scaled}, and $\gamma$ is the Gaussian contour.  This, in turn, implies that 
\[ \bigl| \uPi_{X^n}(H) - \uGamma\bigl\{ J_{X^n}^{1/2}(H - \hat\theta_{X^n}) \bigr\} \bigr| \to 0, \quad \text{in $\prob_\Theta$-probability as $n \to \infty$}, \]
for all compact $H$ with $H \not\ni \Theta$. 
\end{theorem}

The proof of Theorem~\ref{thm:bvn} is lengthy and will be presented elsewhere.  While it follows immediately from \citet{wilks1938} that $\pi_{X^n}(\Theta) \to \gamma(Z)$ in distribution under $\prob_\Theta$, for $Z \sim \nm_d(0,I)$, getting the required in-probability convergence uniformly over parameter values requires significant care.  

Why does $K$ in \eqref{eq:contour.limit} exclude the origin?  It's because there's a sort of singularity, i.e., $\check\pi_{X^n}(0) = \pi_{X^n}(\hat\theta_{X^n})$ is {\em identically 1}, no sampling variability!  Since we know that $\check\pi_{X^n}(0) = \gamma(0)$ exactly, it suffices to focus on $z$ values away from the origin (and from infinity). Similarly, open subsets $H$ that contain $\Theta$ will also eventually contain $\hat\theta_{X^n}$ and, consequently, $\uPi_{X^n}(H)=1$ with high probability when $n$ is large.  Therefore, it suffices to focus on compact $H$ that exclude $\Theta$.

The take-away message is that there's now a mathematically rigorous sense in which this possibilistic IM---despite its imprecision and exact validity---is also statistically efficient.  This also shows that the connection between the possibilistic IMs and Bayes/fiducial solutions highlighted in \citet{martin.isipta2023} for a special class of models holds much more broadly, asymptotically.  Finally, other asymptotically equivalent scaling schemes can be employed in \eqref{eq:im.scaled}.  For example, replacing $J^{-1/2}$ with $n^{-1/2}$ leads to a limiting Gaussian possibility with underlying covariance matrix matching the Cram\'er--Rao lower bound.

\section{Numerical illustrations}
\label{S:num_illustration}

\begin{figure}[t]
\begin{center}
\subfigure[Example~1]{\scalebox{0.44}{\includegraphics{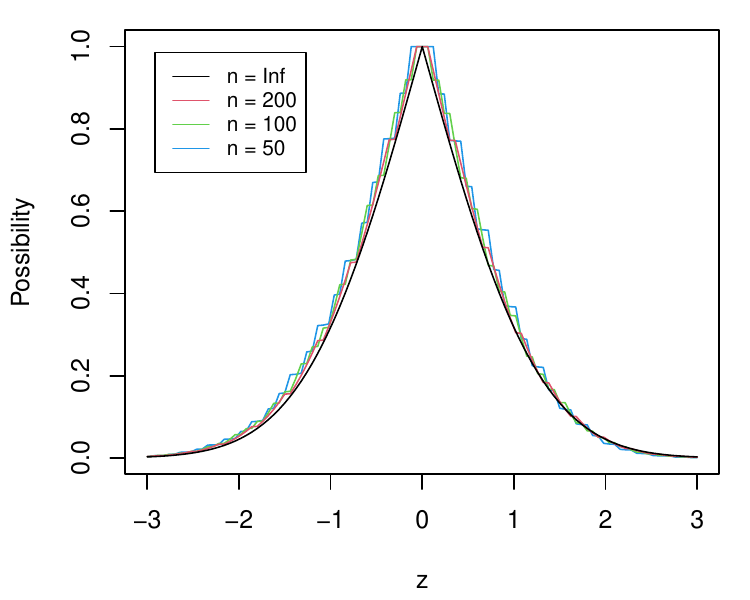}}}
\subfigure[Example~2]{\scalebox{0.44}{\includegraphics{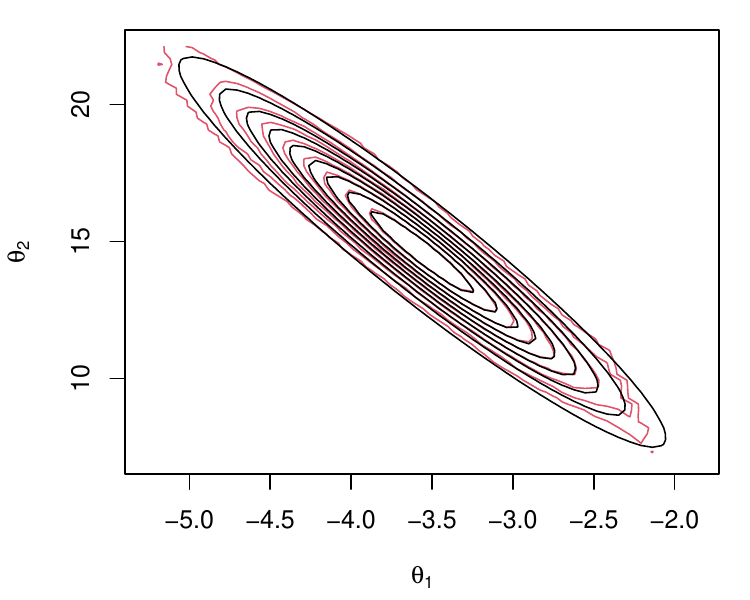}}}
\end{center}
\caption{Possibility contours for two Bernoulli models.}
\label{fig:bernoulli}
\end{figure}

\begin{example} 
{\em iid Bernoulli}.  Figure~\ref{fig:bernoulli}(a) shows the exact centered and scaled possibilistic IM contour for $\Theta$ based on iid Bernoulli data at three different sample sizes, along with the limiting standard Gaussian possibility contour. 
\end{example}

\begin{example}
{\em Logistic regression}.  Consider the data in Table 8.4 of \citet[][p.~252]{ghosh-etal-book} on the relationship between exposure to chloracetic acid and the death of mice. 
A total of $n = 120$ mice are exposed, ten at each of the twelve dose levels, and a binary death indicator is measured.  A simple logistic regression model is fit and Figure~\ref{fig:bernoulli}(b) shows the exact IM possibility contour (red)---which requires expensive but still noisy Monte Carlo evaluations---along with the simple, analytically tractable Gaussian approximation (black).
\end{example}


\section{Conclusion}
\label{S:discuss}

This paper establishes a possibilistic Bernstein--von Mises theorem for the possibilistic IM solution which, in addition to providing new insights and connections to the more familiar Bayesian and fiducial solutions, leads to theoretically and practically important conclusions concerning the IM's efficiency.  Follow-up work will consider extensions of our present result to cases that involve nuisance parameters and/or partial prior information about the unknowns. 



\begin{credits}
\subsubsection{\ackname} RM's research is supported by the U.S.~National Science Foundation, SES--2051225.

\subsubsection{\discintname}
The authors have no competing interests to declare that are
relevant to the content of this article.
\end{credits}

\bibliographystyle{apalike}
\bibliography{mybib}

\begin{thebibliography}{}

\bibitem[Augustin et~al., 2014]{augustin.etal.bookchapter}
Augustin, T., Walter, G., and Coolen, F. P.~A. (2014).
\newblock Statistical inference.
\newblock In {\em Introduction to {I}mprecise {P}robabilities}, Wiley Ser.
  Probab. Stat., pages 135--189. Wiley, Chichester.

\bibitem[Balch et~al., 2019]{balch.martin.ferson.2017}
Balch, M.~S., Martin, R., and Ferson, S. (2019).
\newblock Satellite conjunction analysis and the false confidence theorem.
\newblock {\em Proc. Royal Soc. A}, 475(2227):2018.0565.

\bibitem[Berger, 1984]{berger1984}
Berger, J.~O. (1984).
\newblock The robust {B}ayesian viewpoint.
\newblock In {\em Robustness of {B}ayesian {A}nalyses}, volume~4 of {\em Stud.
  Bayesian Econometrics}, pages 63--144. North-Holland, Amsterdam.
\newblock With comments and with a reply by the author.

\bibitem[Berger et~al., 2009]{bergerbernardosun2009}
Berger, J.~O., Bernardo, J.~M., and Sun, D. (2009).
\newblock The formal definition of reference priors.
\newblock {\em Ann. Statist.}, 37(2):905--938.

\bibitem[Dempster, 1967]{dempster1967}
Dempster, A.~P. (1967).
\newblock Upper and lower probabilities induced by a multivalued mapping.
\newblock {\em Ann. Math. Statist.}, 38:325--339.

\bibitem[Dempster, 2008]{dempster2008}
Dempster, A.~P. (2008).
\newblock The {D}empster--{S}hafer calculus for statisticians.
\newblock {\em Internat. J. Approx. Reason.}, 48(2):365--377.

\bibitem[Den{\oe}ux, 2006]{denoeux2006}
Den{\oe}ux, T. (2006).
\newblock Constructing belief functions from sample data using multinomial
  confidence regions.
\newblock {\em Internat. J. of Approx. Reason.}, 42(3):228--252.

\bibitem[Den{\oe}ux, 2014]{denoeux2014}
Den{\oe}ux, T. (2014).
\newblock Likelihood-based belief function: justification and some extensions
  to low-quality data.
\newblock {\em Internat. J. Approx. Reason.}, 55(7):1535--1547.

\bibitem[Dubois, 2006]{dubois2006}
Dubois, D. (2006).
\newblock Possibility theory and statistical reasoning.
\newblock {\em Comput. Statist. Data Anal.}, 51(1):47--69.

\bibitem[Dubois et~al., 2004]{dubois.etal.2004}
Dubois, D., Foulloy, L., Mauris, G., and Prade, H. (2004).
\newblock Probability-possibility transformations, triangular fuzzy sets, and
  probabilistic inequalities.
\newblock {\em Reliab. Comput.}, 10(4):273--297.

\bibitem[Dubois and Prade, 1988]{dubois.prade.book}
Dubois, D. and Prade, H. (1988).
\newblock {\em Possibility {T}heory}.
\newblock Plenum Press, New York.

\bibitem[Efron, 2013]{efron.cd.discuss}
Efron, B. (2013).
\newblock Discussion: ``{C}onfidence distribution, the frequentist distribution
  estimator of a parameter: a review'' [mr3047496].
\newblock {\em Int. Stat. Rev.}, 81(1):41--42.

\bibitem[Fisher, 1935]{fisher1935a}
Fisher, R.~A. (1935).
\newblock The fiducial argument in statistical inference.
\newblock {\em Ann. Eugenics}, 6:391--398.

\bibitem[Fraser, 1968]{fraser1968}
Fraser, D. A.~S. (1968).
\newblock {\em The {S}tructure of {I}nference}.
\newblock John Wiley \& Sons Inc., New York.

\bibitem[Ghosh et~al., 2006]{ghosh-etal-book}
Ghosh, J.~K., Delampady, M., and Samanta, T. (2006).
\newblock {\em An {I}ntroduction to {B}ayesian {A}nalysis}.
\newblock Springer, New York.

\bibitem[Hannig et~al., 2016]{hannig.review}
Hannig, J., Iyer, H., Lai, R. C.~S., and Lee, T. C.~M. (2016).
\newblock Generalized fiducial inference: a review and new results.
\newblock {\em J. Amer. Statist. Assoc.}, 111(515):1346--1361.

\bibitem[Hose, 2022]{hose2022thesis}
Hose, D. (2022).
\newblock {\em Possibilistic {R}easoning with {I}mprecise {P}robabilities:
  {S}tatistical {I}nference and {D}ynamic {F}iltering}.
\newblock PhD thesis, University of Stuttgart.

\bibitem[Jeffreys, 1946]{jeffreys1946}
Jeffreys, H. (1946).
\newblock An invariant form for the prior probability in estimation problems.
\newblock {\em Proc. Roy. Soc. London Ser. A}, 186:453--461.

\bibitem[Martin, 2015]{plausfn}
Martin, R. (2015).
\newblock Plausibility functions and exact frequentist inference.
\newblock {\em J. Amer. Statist. Assoc.}, 110(512):1552--1561.

\bibitem[Martin, 2018]{gim}
Martin, R. (2018).
\newblock On an inferential model construction using generalized associations.
\newblock {\em J. Statist. Plann. Inference}, 195:105--115.

\bibitem[Martin, 2019]{martin.nonadditive}
Martin, R. (2019).
\newblock False confidence, non-additive beliefs, and valid statistical
  inference.
\newblock {\em Internat. J. Approx. Reason.}, 113:39--73.

\bibitem[Martin, 2021]{imchar}
Martin, R. (2021).
\newblock An imprecise-probabilistic characterization of frequentist
  statistical inference.
\newblock {\tt arXiv:2112.10904}.

\bibitem[Martin, 2022a]{martin.partial}
Martin, R. (2022a).
\newblock Valid and efficient imprecise-probabilistic inference with partial
  priors, {I}. {F}irst results.
\newblock {\tt arXiv:2203.06703}.

\bibitem[Martin, 2022b]{martin.partial2}
Martin, R. (2022b).
\newblock Valid and efficient imprecise-probabilistic inference with partial
  priors, {II}. {G}eneral framework.
\newblock {\tt arXiv:2211.14567}.

\bibitem[Martin, 2023]{martin.isipta2023}
Martin, R. (2023).
\newblock Fiducial inference viewed through a possibility-theoretic inferential
  model lens.
\newblock In Miranda, E., Montes, I., Quaeghebeur, E., and Vantaggi, B.,
  editors, {\em Proceedings of the Thirteenth International Symposium on
  Imprecise Probability: Theories and Applications}, volume 215 of {\em
  Proceedings of Machine Learning Research}, pages 299--310. PMLR.

\bibitem[Martin and Liu, 2013]{imbasics}
Martin, R. and Liu, C. (2013).
\newblock Inferential models: a framework for prior-free posterior
  probabilistic inference.
\newblock {\em J. Amer. Statist. Assoc.}, 108(501):301--313.

\bibitem[Martin and Liu, 2015]{imbook}
Martin, R. and Liu, C. (2015).
\newblock {\em Inferential {M}odels}, volume 147 of {\em Monographs on
  Statistics and Applied Probability}.
\newblock CRC Press, Boca Raton, FL.

\bibitem[Schervish, 1995]{schervish1995}
Schervish, M.~J. (1995).
\newblock {\em Theory of {S}tatistics}.
\newblock Springer-Verlag, New York.

\bibitem[Shafer, 1976]{shafer1976}
Shafer, G. (1976).
\newblock {\em A {M}athematical {T}heory of {E}vidence}.
\newblock Princeton University Press, Princeton, N.J.

\bibitem[Shafer, 1982]{shafer1982}
Shafer, G. (1982).
\newblock Belief functions and parametric models.
\newblock {\em J. Roy. Statist. Soc. Ser. B}, 44(3):322--352.
\newblock With discussion.

\bibitem[van~der Vaart, 1998]{vaart1998}
van~der Vaart, A.~W. (1998).
\newblock {\em Asymptotic {S}tatistics}.
\newblock Cambridge University Press, Cambridge.

\bibitem[Walley, 1991]{walley1991}
Walley, P. (1991).
\newblock {\em Statistical {R}easoning with {I}mprecise {P}robabilities},
  volume~42 of {\em Monographs on Statistics and Applied Probability}.
\newblock Chapman \& Hall Ltd., London.

\bibitem[Wasserman, 1990]{wasserman1990b}
Wasserman, L.~A. (1990).
\newblock Belief functions and statistical inference.
\newblock {\em Canad. J. Statist.}, 18(3):183--196.

\bibitem[Wilks, 1938]{wilks1938}
Wilks, S.~S. (1938).
\newblock The large-sample distribution of the likelihood ratio for testing
  composite hypotheses.
\newblock {\em Ann. Math. Statist}, 9:60--62.

\bibitem[Zabell, 1992]{zabell1992}
Zabell, S.~L. (1992).
\newblock R. {A}. {F}isher and the fiducial argument.
\newblock {\em Statist. Sci.}, 7(3):369--387.

\bibitem[Zadeh, 1978]{zadeh1978}
Zadeh, L.~A. (1978).
\newblock Fuzzy sets as a basis for a theory of possibility.
\newblock {\em Fuzzy Sets and Systems}, 1(1):3--28.

\end{thebibliography}

\end{document}